\documentclass[10pt,reqno]{amsart}
\usepackage[T1]{fontenc}
\usepackage[utf8]{inputenc}
\usepackage{lmodern}
\usepackage{amsmath,amssymb,mathtools,xcolor}
\usepackage[expansion=false]{microtype}
\usepackage[margin=1.5in]{geometry}
\usepackage[colorlinks=true,linkcolor=blue!40!black,
citecolor=blue!40!black,urlcolor=blue!40!black]{hyperref}

\numberwithin{equation}{section}
\newtheorem{theorem}{Theorem}[section]
\newtheorem{lemma}[theorem]{Lemma}
\newtheorem{proposition}[theorem]{Proposition}
\theoremstyle{remark}
\newtheorem{remark}[theorem]{Remark}

\DeclareMathOperator{\diam}{diam}
\DeclareMathOperator{\supp}{supp}
\DeclareMathOperator{\diver}{div}
\newcommand{\R}{\mathbb R}
\newcommand{\eps}{\varepsilon}
\newcommand{\Hone}{H_0^1}
\newcommand{\ip}[2]{\langle #1,#2\rangle}
\newcommand{\norm}[1]{\lVert #1\rVert}

\allowdisplaybreaks[1]
\title[Nodal set of the second eigenfunction]{On the nodal set of the second eigenfunction of the Laplacian in convex domains in $\R^n$}

\author{Jin Sun}
\address{School of Mathematical Sciences, Fudan University, Shanghai, 200433, China}
\email{jsun22@m.fudan.edu.cn}

\date{}

\begin{document}

\begin{abstract}
Let $\Omega$ be a bounded convex domain in $\R^n$ with $n\ge2$.
We prove that the closure of the interior nodal set of any second Dirichlet eigenfunction of $\Omega$ intersects $\partial\Omega$.
The proof combines the ground-state derivative and the quantitative log-concavity estimate of Andrews and Clutterbuck with a variational argument.
\end{abstract}
\maketitle

\section{Introduction}\label{sec:introduction}

Let $\Omega\subset\R^n$ be a bounded connected domain with $n\ge2$.
We consider the Dirichlet eigenvalue problem
\begin{equation}\label{eq:eigenproblem}
\begin{cases}
-\Delta u=\lambda u&\text{in }\Omega,\\
u=0&\text{on }\partial\Omega.
\end{cases}
\end{equation}
It has a discrete spectrum, given by
\[
0<\lambda_1(\Omega)<\lambda_2(\Omega)\le\lambda_3(\Omega)\le\cdots.
\]
For a real eigenfunction $u$, let
\[
Z(u)=\{x\in\Omega:u(x)=0\}
\]
be its interior nodal set.
We say that $u$ has \emph{Payne's property} if $\overline{Z(u)}\cap\partial\Omega\ne\varnothing$.

In 1967, Payne \cite{Payne1967} asked whether the nodal line of a second Dirichlet eigenfunction of a bounded planar domain can be closed.
In \cite{Payne1973}, he proved that the nodal line meets the boundary if the domain is symmetric about a line and convex in the direction perpendicular to it.
Lin \cite{Lin1987} and P\"utter \cite{Putter1990} proved that the nodal line meets the boundary under other symmetry and convexity assumptions on the planar domain.
Jerison \cite{Jerison1991} announced a proof for sufficiently long, thin convex planar domains.
Later, Melas \cite{Melas1992} established Payne's property for smooth convex planar domains, and Alessandrini \cite{Alessandrini1994} extended it to arbitrary bounded convex planar domains.
Jerison \cite{JerisonDiameter1995} estimated the diameter of the nodal line, and Grieser and Jerison \cite{GrieserJerison1996} showed that it approaches a straight line as the eccentricity of the domain tends to infinity.

Much less is known in higher dimensions.
In \cite{Yau1993}, Yau asked about the higher-dimensional analogue of Payne's conjecture and its extension to manifolds with boundary.
Jerison \cite{Jerison1995} proved Payne's property for convex domains whose diameter is sufficiently large relative to their inradius, and Damascelli \cite{Damascelli2000} obtained it under directional convexity and symmetry assumptions.
Freitas and Krej\v{c}i\v{r}\'{i}k \cite{FreitasKrejcirik2008} proved the property for sufficiently thin tubes around curves in Euclidean space.
See also \cite{KrejcirikTusek2015} for thin tubular neighborhoods of hypersurfaces.
The general convex case in dimensions $n\ge3$ is listed as an open problem in \cite[Section~6.3]{Krejcirik2026}.

Without convexity, Payne's property may fail.
Hoffmann-Ostenhof, Hoffmann-Ostenhof and Nadirashvili \cite{Hoffmann1997} constructed bounded multiply connected planar domains with a closed second nodal line.
Dahne, G\'omez-Serrano and Hou \cite{DahneGomezSerranoHou2021} gave a computer-assisted proof for a domain with six holes, and Freitas and Leylekian \cite{FreitasLeylekian2025} later constructed examples with only one hole.
Freitas and Krej\v{c}i\v{r}\'{i}k \cite{FreitasKrejcirik2007} constructed an unbounded simply connected planar domain whose second nodal line is a straight line disjoint from the boundary.
For a closed nodal line on a disk with a rotationally invariant Riemannian metric, see \cite{Freitas2002}.
In higher dimensions, Fournais \cite{Fournais2001} constructed domains with a closed second nodal surface.
Kennedy \cite{Kennedy2013} later constructed examples in every dimension $n\ge3$ that are homeomorphic to a ball and have analytic boundary.
Mukherjee and Saha \cite{MukherjeeSaha2025} used domain perturbations to construct one-parameter families of counterexamples in dimensions $n\ge3$.
Thus, simple connectivity is not sufficient when $n\ge3$.

In this paper, we establish Payne's property for every bounded convex domain in $\R^n$ with $n\ge2$. The main result of this paper is the following.

\begin{theorem}\label{thm:main}
Let $\Omega\subset\R^n$ be a nonempty bounded open convex set, with $n\ge2$.
If $0\ne u\in\Hone(\Omega)$ is a real weak solution of
\[
-\Delta u=\lambda_2(\Omega)u,
\]
then
\[
\overline{Z(u)}\cap\partial\Omega\ne\varnothing.
\]
\end{theorem}

Theorem~\ref{thm:main} therefore resolves the general convex case of Open Problem~5 in \cite[Section~6.3]{Krejcirik2026} for $n\ge3$.

\begin{remark}
Theorem~\ref{thm:main} applies to every nonzero second eigenfunction.
We do not assume that $\lambda_2(\Omega)$ is simple or that $\partial\Omega$ is smooth.
\end{remark}

Our proof is variational.
Let $\phi>0$ be the first Dirichlet eigenfunction of $\Omega$, set $\lambda=\lambda_2(\Omega)$, and consider the shifted energy
\[
Q_\lambda(h)=\int_\Omega|\nabla h|^2-\lambda\int_\Omega h^2.
\]
The variational characterization of $\lambda_2(\Omega)$ gives $Q_\lambda(h)\ge0$ for $h\in\Hone(\Omega)$ orthogonal to $\phi$.
We use the test functions
\begin{equation}\label{eq:A-definition}
A_i u=\phi\,\partial_i(u/\phi)
=\partial_i u-u\,\partial_i\log\phi,
\end{equation}
and write $Au=(A_1u,\ldots,A_nu)$, which we call \emph{ground-state derivative} of $u$.
When $\partial\Omega$ is smooth, by the Hopf lemma, $u/\phi$ is smooth up to the boundary and $A_i u$ vanishes on $\partial\Omega$.
Brascamp and Lieb \cite{BrascampLieb1976} proved that the first Dirichlet eigenfunction of a convex domain is log-concave.
Using the quantitative log-concavity estimate of Andrews and Clutterbuck \cite{AndrewsClutterbuck2011}, we obtain
\begin{equation}\label{eq:intro-energy}
\sum_{i=1}^n Q_\lambda(A_i u)
\le-\frac{2\pi^2}{(\diam\Omega)^2}\norm{Au}_2^2<0.
\end{equation}

Suppose that $Z(u)$ stays away from $\partial\Omega$.
A sufficiently thin shell near the boundary of a convex domain is connected and $u$ has a fixed sign there.
Replacing $u$ by $-u$ if necessary, we may assume that $u$ is positive near the boundary, so that its negative part $w=u^-$ is nonzero and compactly supported in $\Omega$.
The weak equation gives $Q_\lambda(w)=0$ and $Q_\lambda(A_i u,w)=0$.
Subtracting a suitable multiple of $w$ from $A_i u$, we get a function orthogonal to $\phi$ with the same shifted energy as $A_i u$.
Then the variational principle gives $Q_\lambda(A_i u)\ge0$ for every $i$, which is a contradiction to \eqref{eq:intro-energy}.

For a nonsmooth convex domain, we solve $-\Delta U_j=\lambda u$ with zero boundary data on smooth convex domains approximating $\Omega$ from inside.
Passing to the limit in the energy estimates gives $A_i u\in\Hone(\Omega)$ and \eqref{eq:intro-energy}.
The remaining part is similar to the smooth case.

This paper is organized as follows.
Section~\ref{sec:preliminaries} contains the quantitative log-concavity estimate and the construction of smooth convex approximations.
In Section~\ref{sec:derivatives}, we prove two energy identities on smooth domains and pass to an arbitrary convex domain.
We prove Theorem~\ref{thm:main} in Section~\ref{sec:proof}.

\bigskip

\textbf{AI disclosure.}
The author used ChatGPT-6 Pro to explore a variational proof by contradiction based on $Q_\lambda$.
The model suggested the ground-state derivative as a test function and provided a preliminary proof.
The author subsequently introduced Lemmas~\ref{lem:concavity} and~\ref{lem:approximation}, substantially revised and simplified the argument, and rewrote the manuscript.
The author takes full responsibility for the final content of the paper.

\section{Preliminaries}\label{sec:preliminaries}

For an open set $G$, we write $\Hone(G)$ for the closure of $C_c^\infty(G)$ in $H^1(G)$ and understand Dirichlet boundary conditions in this sense.
For $\alpha\in\R$, let
\[
Q_{\alpha,G}(h,k):=\int_G\nabla h\cdot\nabla k-\alpha\int_Ghk,
\qquad Q_{\alpha,G}(h):=Q_{\alpha,G}(h,h).
\]

Consider that $G$ is bounded and connected, and let $\phi>0$ be its first Dirichlet eigenfunction, normalized by $\norm\phi_{L^2(G)}=1$.
The variational characterization of $\lambda_2=\lambda_2(G)$ gives
\begin{equation}\label{eq:minmax}
Q_{\lambda_2,G}(v)\ge0
\qquad\text{for }v\in\Hone(G)\text{ with }\ip v\phi=0.
\end{equation}

The following lemma is a direct consequence of the quantitative log-concavity estimate of Andrews and Clutterbuck \cite[Theorem~1.5]{AndrewsClutterbuck2011}.

\begin{lemma}[Andrews--Clutterbuck]\label{lem:concavity}
Let $G$ be a bounded convex domain, let $d=\diam G$, and let $\phi_G>0$ be its first Dirichlet eigenfunction.
Then
\begin{equation}\label{eq:concavity}
D^2\log\phi_G(x)[\xi,\xi]
\le-\frac{\pi^2}{d^2}|\xi|^2
\qquad(x\in G,\ \xi\in\R^n).
\end{equation}
\end{lemma}

The following standard approximation lemma can be obtained from smooth convex approximation, Dirichlet spectral stability \cite{Daners2003}, and interior elliptic regularity.
For the sake of completeness, we provide the proof here.

\begin{lemma}\label{lem:approximation}
Let $\Omega$ be a bounded convex domain.
There are bounded smooth convex domains $G_j\Subset\Omega$ such that
every compact subset of $\Omega$ is contained in $G_j$ for all
sufficiently large $j$.
Let $\mu_j=\lambda_1(G_j)$, $\mu=\lambda_1(\Omega)$, and let
$\phi_j>0$ and $\phi>0$ be the corresponding $L^2$-normalized ground states.
Then
\begin{equation}\label{eq:ground-convergence}
\mu_j\longrightarrow\mu,
\qquad \phi_j\longrightarrow\phi
\quad\text{in }\Hone(\Omega)\text{ and in }C^m(K)
\end{equation}
for every $K\Subset\Omega$ and every nonnegative integer $m$,
where $\phi_j$ is extended by zero for the $\Hone(\Omega)$ convergence.
\end{lemma}

\begin{proof}
After a translation, we may assume that $0\in\Omega$ and choose
$B_r(0)\subset\Omega\subset B_R(0)$.
The gauge function
\begin{equation}\label{eq:gauge}
p(x):=\inf\{s>0:x\in s\Omega\}
\end{equation}
is convex and positively homogeneous, with $\Omega=\{p<1\}$.
The ball inclusions and subadditivity yield
\[
\frac{|x|}{R}\le p(x)\le\frac{|x|}{r},
\qquad
|p(x)-p(y)|\le\frac{|x-y|}{r}.
\]

Let $\eta\in C_c^\infty(B_1)$ be even and nonnegative with
$\int\eta=1$, and define
\[
p_\eps(x)=\int p(x-\eps y)\eta(y)\,dy.
\]
Then the function $p_\eps$ is smooth and convex.
Jensen's inequality and the Lipschitz bound imply
\[
p\le p_\eps\le p+\frac{\eps}{r}.
\]
Consider 
\[
G_\eps=\{p_\eps<1-\eps\}.
\]
For sufficiently small $\eps>0$, we have
$p_\eps(0)<1-\eps$. Thus, $G_\eps$ is nonempty, and
\[
\overline{G_\eps}\subset\{p\le1-\eps\}\Subset\Omega.
\]
Since $p_\eps$ is smooth and convex, it follows that $G_\eps$ is a bounded smooth convex domain.
Moreover, if $K\Subset\Omega$, then $\max_Kp<1$, which implies that $K\subset G_\eps$ for all sufficiently small $\eps$.
Choosing sufficiently small $\eps_j\downarrow0$ and setting
$G_j=G_{\eps_j}$ gives the required domains.

Extension by zero and the Rayleigh principle give $\mu_j\ge\mu$.
Conversely, approximate $\phi$ in $\Hone(\Omega)$ by functions in
$C_c^\infty(\Omega)$.
Each such function is supported in $G_j$ for all sufficiently large
$j$ and is admissible in the Rayleigh quotient on $G_j$.
Therefore, 
$\limsup_j\mu_j\le\mu$, which implies $\mu_j\to\mu$.

Since the zero extensions of $\phi_j$ are bounded in $\Hone(\Omega)$,
every subsequence has a further subsequence converging
weakly in $\Hone(\Omega)$ and strongly in $L^2(\Omega)$ to some
$\phi_*\ge0$ with $\|\phi_*\|_{L^2(\Omega)}=1$.
The Rayleigh principle and weak lower semicontinuity give
\[
\mu\le\|\nabla\phi_*\|_{L^2(\Omega)}^2
\le\lim_{j\to\infty}\mu_j=\mu.
\]
Thus, $\phi_*$ is a nonnegative normalized ground state, which means $\phi_*=\phi$.
Consequently, the whole sequence converges weakly in $\Hone(\Omega)$
and strongly in $L^2(\Omega)$ to $\phi$.
Since
\[
\|\nabla\phi_j\|_{L^2(\Omega)}^2
=\mu_j\longrightarrow\mu
=\|\nabla\phi\|_{L^2(\Omega)}^2,
\]
the convergence in $\Hone(\Omega)$ is strong.

Finally, fix $K\Subset U\Subset\Omega$ with $U$ open.
For all sufficiently large $j$, the equations
$-\Delta\phi_j=\mu_j\phi_j$ hold on $U$.
Therefore, interior elliptic estimates and Sobolev embedding yield $\phi_j\to\phi$ in $C^m(K)$ for every nonnegative integer $m$.
\end{proof}

\section{Ground-state derivatives}\label{sec:derivatives}

We write $\partial_i=\partial/\partial x_i$ and
$\partial_{ik}=\partial_i\partial_k$ for derivatives.
In the following, we consider the ground-state derivative $A$ in \eqref{eq:A-definition}, given by 
\[
A u=\phi\,\nabla(u/\phi)
=\nabla u-u\,\nabla\log\phi.
\]

First, we give two energy identities for smooth functions vanishing on the boundary.
The second identity will allow us to approximate $u$ by solutions of a Poisson equation on smooth domains.

\begin{lemma}\label{lem:identities}
Let $G$ be a bounded connected domain with $C^\infty$ boundary, let $\phi>0$ be its first Dirichlet eigenfunction, and set $\psi=\log\phi$.
If $h\in C^\infty(\overline G)$ and $h=0$ on $\partial G$, then $A_i h\in C^\infty(\overline G)\cap\Hone(G)$ and
\begin{align}
\norm{Ah}_{L^2(G)}^2
&=\norm{\nabla h}_{L^2(G)}^2-\lambda_1(G)\norm h_{L^2(G)}^2,
\label{eq:factorization}\\
\sum_{i=1}^n\norm{\nabla A_i h}_{L^2(G)}^2
&=\norm{\Delta h}_{L^2(G)}^2-\lambda_1(G)\norm{\nabla h}_{L^2(G)}^2
+2\int_GD^2\psi[Ah,Ah].
\label{eq:second-identity}
\end{align}
\end{lemma}

\begin{proof}
By the Hopf lemma, we obtain that $f=h/\phi$ is smooth on $\overline G$.
It follows that $A_i h=\phi\,\partial_i f$ is smooth on $\overline G$ and vanishes on $\partial G$.
Since
\[
|\nabla h|^2
=\phi^2|\nabla f|^2+\nabla\phi\cdot\nabla(\phi f^2),
\]
integration by parts yields 
\begin{align*}
\norm{\nabla h}_{L^2(G)}^2
&=\int_G\phi^2|\nabla f|^2+\int_G\nabla\phi\cdot\nabla(\phi f^2)
\\
&=\norm{Ah}_{L^2(G)}^2+\lambda_1(G)\int_G\phi^2 f^2
\\
&=\norm{Ah}_{L^2(G)}^2+\lambda_1(G)\norm h_{L^2(G)}^2,
\end{align*}
which proves \eqref{eq:factorization}.

Note that the equation for $\phi$ gives
\[
\Delta\psi+|\nabla\psi|^2=-\lambda_1(G).
\]
Differentiating this identity, we obtain
\[
\Delta(\partial_i\psi)=-2\sum_k(\partial_k\psi)(\partial_{ik}\psi).
\]
Therefore,
\begin{equation}\label{eq:commutator}
\begin{aligned}
-\Delta(A_i h)
&=-A_i(\Delta h)+h\,\Delta(\partial_i\psi)
+2\sum_k(\partial_k h)(\partial_{ik}\psi)\\
&=-A_i(\Delta h)+2\sum_k(\partial_{ik}\psi)A_k h.
\end{aligned}
\end{equation}
For the term $-A_i(\Delta h)$, since $A h=0$ on $\partial G$ and
\[
\diver (A h)+\nabla\psi\cdot Ah=\Delta h+\lambda_1(G) h,
\]
integration by parts yields
\begin{align*}
-\int_GA h\cdot A(\Delta h)
&=\int_G\Delta h\,\bigl(\diver A h+\nabla\psi\cdot A h\bigr)\\
&=\int_G\Delta h\,(\Delta h+\lambda_1(G) h)\\
&=\norm{\Delta h}_{L^2(G)}^2-\lambda_1(G)\norm{\nabla h}_{L^2(G)}^2.
\end{align*}
Multiplying \eqref{eq:commutator} by $A_i h$, integrating over $G$ and summing over $i$, we obtain \eqref{eq:second-identity}.
\end{proof}
\begin{remark}
Although the coefficients in \eqref{eq:commutator} may be singular at $\partial G$, the products used above extend smoothly to $\overline G$.
Indeed,
\[
-A_i h\,A_i(\Delta h)
=-\phi(\partial_i f)\,\partial_i(\Delta h)
+(\partial_i\phi)(\partial_i f)\,\Delta h
\]
and
\[
(\partial_{ik}\psi)(A_i h)(A_k h)
=\bigl(\phi\,\partial_{ik}\phi
-(\partial_i\phi)(\partial_k\phi)\bigr)
(\partial_i f)(\partial_k f).
\]
Therefore, integration by parts is valid.
\end{remark}

If $G$ is also convex and $d=\diam G$, then Lemma~\ref{lem:concavity} and \eqref{eq:second-identity} imply
\begin{equation}\label{eq:smooth-estimate}
\sum_{i=1}^n\norm{\nabla A_i h}_{L^2(G)}^2
+\frac{2\pi^2}{d^2}\norm{Ah}_{L^2(G)}^2
\le\norm{\Delta h}_{L^2(G)}^2-\lambda_1(G)\norm{\nabla h}_{L^2(G)}^2.
\end{equation}

We next prove that $A_i u\in\Hone(\Omega)$ and establish the energy estimate for any bounded convex domain without boundary regularity assumption.

\begin{proposition}\label{prop:energy}
Let $\Omega$ be a bounded convex domain, let $\phi>0$ be its first Dirichlet eigenfunction, and set $\mu=\lambda_1(\Omega)$ and $D=\diam\Omega$.
Suppose that $u\in\Hone(\Omega)$ satisfies $-\Delta u=\lambda u$, where $\lambda>\mu$.
Then $A_i u\in\Hone(\Omega)$ for each $1\le i\le n$,
\begin{equation}\label{eq:admissibility}
\norm{Au}_{L^2(\Omega)}^2=(\lambda-\mu)\norm u_{L^2(\Omega)}^2,
\end{equation}
and
\begin{equation}\label{eq:negative-energy}
\sum_{i=1}^nQ_{\lambda,\Omega}(A_i u)
\le-\frac{2\pi^2}{D^2}\norm{Au}_{L^2(\Omega)}^2.
\end{equation}
\end{proposition}

\begin{proof}
Let $G_j$, $\mu_j$ and $\phi_j$ be as in Lemma~\ref{lem:approximation}.
Since $u$ may not vanish on $\partial G_j$, we first approximate it by functions in $\Hone(G_j)$ to which Lemma~\ref{lem:identities} applies.
Let $U_j\in\Hone(G_j)$ solve
\begin{equation}\label{eq:poisson}
-\Delta U_j=\lambda u\quad\text{in }G_j.
\end{equation}
Existence and uniqueness follow from the Lax--Milgram theorem.
After extending $U_j$ by zero, the weak equations for $u$ and $U_j$ give
\[
\int_\Omega\nabla(u-U_j)\cdot\nabla v=0,
\qquad \forall\,v\in\Hone(G_j).
\]
Thus, $U_j$ is the orthogonal projection of $u$ onto $\Hone(G_j)$ with respect to the Dirichlet inner product on $\Hone(\Omega)$, and
\begin{equation}\label{eq:projection}
\norm{\nabla(u-U_j)}_{L^2(\Omega)}
=\inf_{z\in\Hone(G_j)}\norm{\nabla(u-z)}_{L^2(\Omega)}
\longrightarrow0,
\end{equation}
since $u$ can be approximated in $\Hone(\Omega)$ by functions in $C_c^\infty(\Omega)$. 
Therefore, Poincar\'e's inequality gives $U_j\to u$ in $\Hone(\Omega)$.

Boundary regularity gives $U_j\in C^\infty(\overline{G_j})$.
Moreover, $U_j-u$ is harmonic in $G_j$.
Thus, interior estimates yield
\begin{equation}\label{eq:U-convergence}
U_j\longrightarrow u\quad\text{in }C^m(K),
\qquad \forall\,K\Subset\Omega,\ m\ge0.
\end{equation}

Set $A^{(j)}=\nabla-\nabla\log\phi_j$ and $V^{(j)}=A^{(j)}U_j$.
By Lemma~\ref{lem:identities}, each component of $V^{(j)}$ belongs to $\Hone(G_j)$.
We extend $V^{(j)}$ by zero to $\Omega$.
Using \eqref{eq:ground-convergence}, \eqref{eq:factorization} and \eqref{eq:projection}, we obtain
\begin{equation}\label{eq:V-norm}
\norm{V^{(j)}}_{L^2(\Omega)}^2
=\norm{\nabla U_j}_{L^2(\Omega)}^2-\mu_j\norm{U_j}_{L^2(\Omega)}^2
\longrightarrow(\lambda-\mu)\norm u_{L^2(\Omega)}^2.
\end{equation}
Since $\diam G_j\le D$, applying \eqref{eq:smooth-estimate} to $U_j$ and using \eqref{eq:poisson}, we obtain
\begin{equation}\label{eq:uniform-estimate}
\sum_i\norm{\nabla V^{(j)}_i}_{L^2(\Omega)}^2
+\frac{2\pi^2}{D^2}\norm{V^{(j)}}_{L^2(\Omega)}^2
\le\lambda^2\int_{G_j}u^2-\mu_j\norm{\nabla U_j}_{L^2(\Omega)}^2.
\end{equation}
The right-hand side is bounded above by $\lambda^2\norm u_{L^2(\Omega)}^2$.
Hence, for each $1\le i\le n$, the sequence $\{V^{(j)}_i\}$ is bounded in $\Hone(\Omega)$.
Passing to a common subsequence, we may assume that
\[
V^{(j)}_i\rightharpoonup V_i\quad\text{in }\Hone(\Omega)
\]
and
\[
V^{(j)}_i\longrightarrow V_i\quad\text{in }L^2(\Omega)
\]
for every $i$.
Set $V=(V_1,\ldots,V_n)$.

Fix $K\Subset\Omega$.
Since $\phi$ has a positive minimum on $K$, \eqref{eq:ground-convergence} and \eqref{eq:U-convergence} imply
\[
V^{(j)}=\nabla U_j-U_j\frac{\nabla\phi_j}{\phi_j}
\longrightarrow\nabla u-u\frac{\nabla\phi}{\phi}=Au
\quad\text{uniformly on }K.
\]
Thus, $V=Au$ almost everywhere.
The weak convergence of the components gives $A_i u\in\Hone(\Omega)$ for every $i$, and the strong $L^2$ convergence and \eqref{eq:V-norm} give 
\[
\norm{Au}_{L^2(\Omega)}^2 = \lim_{j\to\infty}\norm{V^{(j)}}_{L^2(\Omega)}^2 = (\lambda-\mu)\norm{u}_{L^2(\Omega)}^2,
\]
which proves \eqref{eq:admissibility}.

Finally, together with \eqref{eq:ground-convergence} and \eqref{eq:projection}, 
\[
\lambda^2\int_{G_j}u^2-\mu_j\norm{\nabla U_j}_{L^2(\Omega)}^2\longrightarrow \lambda(\lambda-\mu)\norm u_{L^2(\Omega)}^2.
\]
Therefore,
\[
\begin{aligned}
\sum_i\norm{\nabla A_i u}_{L^2(\Omega)}^2
+\frac{2\pi^2}{D^2}\norm{Au}_{L^2(\Omega)}^2
&\le\liminf_j\left(
\sum_i\norm{\nabla V^{(j)}_i}_{L^2(\Omega)}^2
+\frac{2\pi^2}{D^2}\norm{V^{(j)}}_{L^2(\Omega)}^2\right)\\
&\le\lambda(\lambda-\mu)\norm u_{L^2(\Omega)}^2
=\lambda\norm{Au}_{L^2(\Omega)}^2,
\end{aligned}
\]
which is equivalent to \eqref{eq:negative-energy}.
\end{proof}

\begin{remark}
The approximation \eqref{eq:poisson} only uses the given eigenfunction $u$, so no simplicity assumption on $\lambda$ or choice of eigenfunctions on $G_j$ is needed.
\end{remark}

\section{\texorpdfstring{Proof of Theorem~\ref{thm:main}}
{Proof of the main theorem}}\label{sec:proof}

We will use a compactly supported nodal part of $u$ to make the test functions orthogonal to the first eigenfunction while preserving their shifted energy.
Write $u^-=\max\{-u,0\}$ for the negative part of $u$.

\begin{lemma}\label{lem:nodal-part}
Let $\Omega$ be a bounded connected domain, let $\lambda=\lambda_2(\Omega)$, and let $u$ be an eigenfunction associated with $\lambda$.
Suppose that $w=u^-$ is nonzero and $\supp w\Subset\Omega$.
If $v=\partial_i u-au\in\Hone(\Omega)$ for some $i$ and for some $a\in C^\infty(\Omega)$, then $Q_\lambda(v)\ge0$.
\end{lemma}

\begin{proof}
The Sobolev chain rule gives $w\in\Hone(\Omega)$ and
\[
\nabla w=-\mathbf1_{\{u<0\}}\nabla u
\quad\text{almost everywhere}.
\]
Since $uw=-w^2$, testing the equation for $u$ against $w$ yields $Q_\lambda(w)=0$.
We also have
\[
u\nabla w=w\nabla u
\quad\text{almost everywhere}.
\]
Let $\phi>0$ be the first eigenfunction.
Since $w\ge0$ and $w\ne0$, we have $\ip w\phi>0$.

We next show that $Q_\lambda(v,w)=0$.
Choose an open set $U\Subset\Omega$ containing $\supp w$.
Note that $w\in\Hone(U)$.
Since $a$ is smooth near $\overline U$, we also have $aw\in\Hone(U)$.
Differentiating the equation for $u$ in the interior and testing against $w$ gives 
\[
Q_\lambda(\partial_i u,w)=\int_U\nabla(\partial_i u)\cdot\nabla w-\lambda\int_U(\partial_i u)w=0.
\]
Testing the original equation against $aw$, we obtain
\[
\begin{aligned}
Q_\lambda(au,w)
=Q_\lambda(u,aw)
+\int_U(u\nabla w-w\nabla u)\cdot\nabla a
=0.
\end{aligned}
\]
Thus, we have $Q_\lambda(v,w)=0$.

Set
\[
c=\frac{\ip v\phi}{\ip w\phi},\qquad z=v-cw.
\]
Then $z\in\Hone(\Omega)$ and $z\perp\phi$.
By \eqref{eq:minmax},
\[
0\le Q_\lambda(z)
=Q_\lambda(v)-2cQ_\lambda(v,w)+c^2Q_\lambda(w)
=Q_\lambda(v),
\]
which completes the proof.
\end{proof}

Now we are ready to prove our main theorem.
\begin{proof}[Proof of Theorem~\ref{thm:main}]
Let $\lambda=\lambda_2(\Omega)$ and let $\phi>0$ be the first eigenfunction.
Suppose for contradiction that $\overline{Z(u)}\cap\partial\Omega=\varnothing$.
Then $\overline{Z(u)}$ is a compact subset of $\Omega$.

After translating an interior point to the origin, let $p$ be the gauge function in \eqref{eq:gauge}.
Choose $t\in(0,1)$ such that $\overline{Z(u)}\subset\{p<t\}$.
Since $n\ge2$, the shell $\{t<p<1\}$ is connected.
The function $u$ is continuous and does not vanish there, so it has a fixed sign.
Replacing $u$ by $-u$ if necessary, we may assume that $u$ is positive on the shell.
Hence
\[
\supp u^-\subset\{p\le t\}\Subset\Omega.
\]
The weak eigenvalue equations give $(\lambda-\lambda_1)\ip u\phi=0$, and thus $u\perp\phi$.
Since $u\ne0$, we have $w=u^-\ne 0$.

By Proposition~\ref{prop:energy}, $A_i u\in\Hone(\Omega)$ for every $i$.
Applying Lemma~\ref{lem:nodal-part} with $a=\partial_i\log\phi$, we obtain $Q_\lambda(A_i u)\ge0$ for every $i$.
On the other hand, \eqref{eq:admissibility} and \eqref{eq:negative-energy} give
\[
0\le\sum_iQ_\lambda(A_i u)
\le-\frac{2\pi^2}{(\diam\Omega)^2}
(\lambda_2-\lambda_1)\norm u_2^2<0.
\]
This contradiction proves the theorem.
\end{proof}

\begin{remark}
The second eigenvalue enters the proof only through \eqref{eq:minmax}.
For higher eigenvalues, orthogonality to the first eigenfunction does not make the shifted energy nonnegative.
\end{remark}

\section*{Acknowledgments}
We thank Prof. Bobo Hua for his guidance and constant support.

\bibliographystyle{amsalpha}
\bibliography{ref}
\end{document}